\documentstyle[12pt]{article}
\def\Bbb#1{{\mathchoice{\mbox{\bf #1}}{\mbox{\bf #1}}%
{\mbox{$\scriptstyle \bf #1$}}{\mbox{$\scriptscriptstyle \bf #1$}}}}

\def\C{\Bbb C}
\def\Z{\Bbb Z}
\def\T{\Bbb T}

\def\E{\Bbb E}
\def\O{{\cal O}}
\def\sgn{{\rm sgn}}

\def\A{{\cal A}}
\def\cT{{\cal T}}
\def\cF{{\cal F}}

\begin{document}
\title{On a weak type $(1,1)$ inequality for 
a maximal conjugate function}
\author{Nakhl\'e H. \ Asmar
 and 
 Stephen J.\ Montgomery-Smith
}
\date{}
\maketitle
\section{Introduction}
\newtheorem{maintheorem}{Theorem}[section]
\newtheorem{remarks}[maintheorem]{Remarks}

Throughout this paper, $N$ denotes a fixed but arbitrary 
positive integer, $\T$ denotes the
circle group, and $\T^N$
denotes the product of $N$ copies of $\T$.
The normalized Lebesgue measure on $\T^N$ will be 
symbolized by $P$.  For a measurable function
$f$, we let $\|f\|^*_1=\sup_{y>0}y \lambda_f(y)$
where
$\lambda_f (y)=P\left(\{x\in \T^N:\ |f(x)|>y\}
\right)$.  The integers will be denoted
by $\Z$ and the complex numbers by $\C$.

Let $\cF_n=\sigma(e^{i\theta_1},e^{i\theta_2},\ldots,e^{i\theta_n})$
denote the $\sigma$-algebra on $\T^N$ generated by the first $n$ coordinate functions.
For $f\in L^1(\T^N)$, the conditional expectation of $f$
with respect to $\cF_n$ will be denoted 
$\E(f|\cF_n)$.  Let
$$d_0(f)=\E(f|\cF_0)=\int_{\T^N} fdP,$$
and for $j=1,\ldots, N$, let
$d_j(f)=\E(f|\cF_j)-\E(f|\cF_{j-1})$.
We have the martingale difference decomposition 
\begin{equation}
f=\sum_{j=0}^N d_j(f).
\label{martingale-difference}
\end{equation}
Consider the  
maximal function corresponding to 
(\ref{martingale-difference})
\begin{equation}
D(f)=
\sup_{1\leq n\leq N}
\left|
\sum_{j=0}^n d_j(f)
\right|
=
\sup_{1\leq n\leq N}
\left|
\E(f|\cF_n)
\right|.
\label{max-martingale}
\end{equation}
A well-known weak type $(1,1)$ maximal inequality due to
Doob states that there is 
a constant $a$ independent of $f$ and $N$ such that
\begin{equation}
\|Df\|^*_1\leq a \|f\|_1.
\label{doob's-inequality}
\end{equation}
Now we recall the conjugate function
operator $f\mapsto \widetilde{f}$,
defined for all $f\in L^2(\T)$ by
the multiplier relation
$$\widehat{\widetilde{f}}(n)=-i\sgn (n)\widehat{f}(n),
\ {\rm for\ all}\ n\in \Z. $$
By Kolmogorov's Theorem \cite[Chap.\ IV, Theorem (3.16)]{zyg},
the operator $f\mapsto \widetilde{f}$ is of weak type
$(1,1)$.  

Denote an element of $\T^N$ by $(\theta_1,\theta_2,\ldots,\theta_N)$.
Let $H_j$ denote the one-dimensional conjugate function
operator defined for functions on $\T^N$ with respect to the
$\theta_j$ variable.  As an operator on $L^2(\T^N)$,
$H_j$ is given by the multiplier relation
$\widehat{H_j(f)}(z_1,z_2,\ldots,z_N)=
-i\sgn(z_j)\widehat{f}(z_1,z_2,\ldots,z_N)$,
for all $(z_1,z_2,\ldots,z_N)\in \Z^N$.
Plainly, the operators $H_j$, $j=1,\ldots,N$, are
of weak type $(1,1)$ on $L^1(\T^N)$ with the same constant
as in Kolmogorov's theorem for $L^1(\T)$.
The conjugate function that we consider is defined 
for all 
$f\in L^1(\T^N)$ by
\begin{equation}
H(f)=\sum_{j=1}^N H_j(d_j(f)).
\label{conj-function}
\end{equation}
Since both $H_j$ and $d_j$ are multipliers, they commute.
We have
\begin{equation}
H(f)=\sum_{j=1}^N d_j(H(f)).
\label{conj-mart-dif}
\end{equation}
The maximal function that we are interested in is defined by
\begin{equation}
M(f)=\sup_{1\leq n\leq N}
\left|
\sum_{j=1}^n d_j(H_j(f))
\right|
=D(H(f)),
\label{max-conj-ft}
\end{equation}
where $D$ is as in (\ref{max-martingale}).  
Thus $M$ is the composition of two operators of weak type
$(1,1)$.  (The fact that $H$ is of weak type $(1,1)$
is known, and will not be needed in the proofs.  See 
Remarks \ref{remarks} (a), below.  This fact will also follow 
from our main theorem.)
Our goal is to
prove the following result.  
\begin{maintheorem}
There is a constant $A$ independent of $N$ such that
for all $f\in L^1(\T^N)$ we have
\begin{equation}
\|Mf\|^*_1\leq A \|f\|_1,
\label{desired-inequality}
\end{equation}
where $M$ is the maximal operator given by 
(\ref{max-conj-ft}).
\label{maintheorem}
\end{maintheorem}
The proof of this theorem is presented in the following section,
and is of independent interest.
We will show that by changing the 
time in the Brownian motion that Burkholder, Gundy, and Silverstein
used in \cite{bgs} from a continuous range $[0,\infty)$
to a semi-continuous range $\{1,2,\ldots\}\times [0,\infty)$,
the proofs in \cite{bgs} can be carried out
on $\T^N$,
yielding inequalities which are independent of $N$
(e.g., the ``good $\lambda$'' inequality).

We end this section with some remarks
concerning the operator $H$ that will not be used in the sequel.
\begin{remarks}
{\rm
(a)  The operator $f\mapsto Hf$ that we defined in 
(\ref{conj-mart-dif})
is a conjugate function operator
 of the kind that was introduced
and studied by Helson \cite{hel}.
Helson's definition is in terms of orders on the
dual group $\Z^N$.  In our case, the operator
$H$ can be recast in terms of a lexicographic order on
$\Z^N$.  As shown in \cite{hel}, the operator
$H$ is bounded from $L^1(\T^N)$ into $L^p(\T^N)$,
for any $0<p<1$.  Indeed it is of weak type $(1,1)$
(see \cite[Theorem 4.3]{ams}).\\
(b)
We proved in \cite[Theorem 5.4]{ams} that the 
square function $Sf=\left(
\sum_{j=1}^N
|H_j(d_j(f))|^2
\right)^{1/2}$
is of weak type $(1,1)$.  
It is known that under certain conditions
on the martingale, the weak type estimates for 
the square function and the maximal function
are equivalent (see, for example, 
\cite{bg}, Assumptions A1-A3).  The martingales that we
are studying do not satisfy these conditions,
and so 
(\ref{desired-inequality})
does not follow from the weak $(1,1)$
estimates for the square function,
by using general facts from probability theory.
}
\label{remarks}
\end{remarks}
\section{Proof of Theorem 1.1}
\newtheorem{lemma5.1}{Lemma}[section]
\newtheorem{lemma5.2}[lemma5.1]{Lemma}
\newtheorem{lemma5.3}[lemma5.1]{Lemma}
\newtheorem{lemma5.4}[lemma5.1]{Lemma}
For clarity's sake, we start with an outline of the proof,
setting in the process our notation, and describing our 
generalization of the methods in \cite{bgs}.

It is enough to prove (\ref{desired-inequality}) with $f\in S(\T^N)$, 
the space of trigonometric polynomials
 on $\T^N$.  We may also assume that
 $f$ is real-valued and that $d_0(f)=0$.  Write
$$ f(\theta_1,\dots,\theta_N) = \sum a_{j_1,\dots,j_N} \theta_1^{j_1}
   \dots \theta_N^{j_N},$$
and extend 
$f$\ to a function on 
$\C^N$\ that
is harmonic in each variable as follows 
$$ f(r_1 \theta_1, \dots, r_N \theta_N) = \sum a_{j_1,\dots,j_N} 
   r_1^{|j_1|} \theta_1^{j_1}
   \dots r_N^{|j_N|} \theta_N^{j_N} $$
where $r_n$\ is a nonnegative real number, 
and $|\theta_n| \in \T = \{z:
|z| = 1\}$.  In this notation, 
the $n$-th term in the martingale difference decomposition 
of $f$ becomes
$$d_n(f)=\sum_{\stackrel{j_1,j_2,\ldots,j_n}{j_n\neq 0}} 
a_{j_1,\ldots ,j_n,0,\ldots,0} \theta_1^{j_1} \ldots \theta_n^{j_n} .$$
Since by assumption $d_0(f)=0$, it follows that
\begin{equation}
d_n (f) (r_1 \theta_1, \ldots , r_{n-1} \theta_{n-1}, 0)=0
\label{value-at-0}
\end{equation}
for all $n=0,1, \ldots, N.$\\
 
The approach that we take is to consider a martingale 
on a time structure
that is part continuous and part discrete.  
Our notion
of time is $\cT = \{1,2,\dots,N\} \times 
[0,\infty[$\ with the order $(m,s) < (n,t)$\
if and only if $m<n$\ or $m =n$\ and $s < t$.  Construct $N$\ independent
complex Brownian motions $c_{n,t} = a_{n,t} + i b_{n,t}$\ 
($1 \le n \le N,\ t\geq 0$)
each one starting at $0$.  
Define stopping times $\tau_n = \inf\{\, t:
|c_{n,t}| \ge 1\}$.  

Define an increasing family of sigma fields 
$(\A_{(n,t)} : (n,t) \in 
\cT )$, where
$ \A_{(n,t)} $\ is the sigma field generated by the functions $c_{m,s}$\
for $(m,s) \le (n,t)$.  
Then we define a process over our new time structure by:
\begin{eqnarray}
 F_{n,t} &=&  f(c_{1,\tau_1},\ldots,c_{n-1,\tau_{n-1}},
 c_{n,\tau_n\wedge t},0,\ldots,0) \nonumber   \\
 &=& \sum_{k=0}^{n-1}
   d_k(f) (c_{1,\tau_1},\dots,c_{k,\tau_k})
   + d_n (f) (c_{1,\tau_1},\dots,c_{n,t \wedge \tau_n}) .
\label{brownian-martingale}
\end{eqnarray}
Since $\tau_n < \infty$\ a.s., it follows that a.s.,
for sufficiently large
$(n,t)$, we have $F_{n,t} = F_\infty$, where
$$ F_\infty = \sum_{k=0}^N
   d_k(c_{1,\tau_1},\dots,c_{k,\tau_k}) = 
   f(c_{1,\tau_1},\dots,c_{N,\tau_N}).$$
We will show that the family of functions
$(F_{n,t})$ is a martingale relative to
$ \A_{(n,t)} $. To be able
to use results from the classical theory of 
martingales,
it is convenient to 
label the family $(F_{n,t})$ by a continuous 
time parameter.  This can be done by forming
an order
preserving bijection between $\cT \cup \{\infty\}$\ and 
$[0,N]$\ as follows:
$$\phi(n,t) = n-1 + t/(t+1),\ {\rm and}\  \phi(\infty) = N.$$
Because $c_{n,t}$\ is a.s.\ continuous in $t$, and 
also $\tau_n < \infty$\ a.s., it follows that
$F_{\phi^{-1}(t)}$ is a
continuous time martingale on $[0,N]$.  
Let $\tilde F_{n,t}$\ be constructed from $H f$\ as in (\ref{brownian-martingale}).
Define the Brownian maximal function
$$ F^* = \sup_{0 \le t \le N} |F_{\phi^{-1}(t)}| ,$$
and let $\tilde F^*$ be defined similarly by using $\tilde F_{n,t}$.  
The proof of the desired inequality (\ref{desired-inequality}) 
will proceed in four steps.\\
Step 1:   $\left\| F_\infty \right\|_1 = \left\| f \right\|_1$;\\
Step 2:   $\left\| F^* \right\|^*_{1,\infty} \le   \left\| F_\infty \right\|_1$;\\
Step 3:   $\left\| \tilde F^* \right\|^*_{1,\infty} \le c \, \left\|
            F^* \right\|^*_{1,\infty}$;\\
Step 4:   $\left\| M f \right\|^*_{1,\infty} \le \left\| \tilde F^* 
            \right\|^*_{1,\infty}$.\\

We now proceed with the proofs.
Suppose that 
$c_t=a_t + i b_t$\ is a complex Brownian motion starting at $0$.  
Let $A_t$\ be the sigma field generated by $c_s$\ for $s \le t$.
Let
$\tau = \inf\{ t: |c_t| \ge 1\}$.\\
Suppose that $v$ is 
 a real-valued trigonometric polynomial on 
$\T=\{|z|=1\}$, and extend $v$ to be 
harmonic on $\C$.  It follows
from \cite[Theorem 4.1]{doob2} that $v(c_t)$
is a martingale, and $v(c_t)$ is $A_t$-measurable.
The following lemma, is a simple consequence of this fact and 
Doob's Optional Stopping Theorem.  
\begin{lemma5.1}
With the above notation, 
if $\mu$\ is a stopping time such that
$\mu \le \tau$, then
$$ \E( v(c_\mu) | A_t) = v(c_{t \wedge \mu}) .$$
\label{lemma5.1}
\end{lemma5.1}
%
%
%
%
%
%
%
%
%
%
%
%
\cite[Theorem (3.2), p.65]{ry}.  We have
%
%
%
%
%
%
%
Using Lemma (\ref{lemma5.1}), we can establish a basic property of the 
functions $(F_{n,t})$.
\begin{lemma5.2}
In the above notation, we have that
$\E(F_\infty | \A_{n,t}) = F_{n,t}$, and hence that
$(F_{n,t})$\ is a martingale. Consequently, $(F_{\phi^{-1}(t)})$ 
is a  continuous time martingale for
$t\in [0,N]$.
\label{lemma5.2}
\end{lemma5.2}
{\bf Proof.}  First, it is clear that if $ k < n$, then
$$ \E(d_k(c_{1,\tau_1},\dots,c_{k,\tau_k}) | \A_{n,t})
   = d_k(c_{1,\tau_1},\dots,c_{k,\tau_k}), $$
 because $d_k(c_{1,\tau_1},\dots,c_{k,\tau_k})$\ is $\A_{n,t}$\
measurable.  Also, if $k > n$, then
$$ \E(d_k(c_{1,\tau_1},\dots,c_{k,\tau_k}) | \A_{n,t})
   =
   \E(\E(d_k(c_{1,\tau_1},\dots,c_{k,\tau_k}) | \A_{k,0}) | \A_{n,t})
   = 0 ,$$
by Lemma (\ref{lemma5.1})
and (\ref{value-at-0}).  
Similarly, by the same lemma, it also follows that if
$k = n$, then
$$ \E(d_k(c_{1,\tau_1},\dots,c_{k,\tau_k}) | \A_{n,t}) 
   =
   d_k(c_{1,\tau_1},\dots,c_{k,t \wedge \tau_k}) $$
and hence $\E(F_\infty | \A_{n,t}) = F_{n,t}$.  This proves that 
$(F_{n,t})$\ is a martingale.  The rest of the lemma is obvious.\\
{\bf Proof of Steps 1, 2, 4}
Because of Lemma (\ref{lemma5.2}), Step 2 follows from 
Doob's Maximal Inequality for continuous time
martingales (see \cite[Chapter VII, Section 11]{doob}).
Step 1 also follows from the uniform distribution
of Brownian motion over $\T$ (see \cite[Corollary 3.6.2]{peter}).
Step 4 is also a consequence of the same property of Brownian motion. 
We give details.  We have
\begin{eqnarray*}
\tilde F^*      & = & \sup_{(n,t)}| \tilde F_{n,t} | \geq
                                        \sup_n | \tilde F_{n,\tau_n} |\\
                & = & \sup_n \left| 
                                \sum_{m=0}^n 
                                 H_m (d_m (f)) 
                                (c_{1,\tau_1},\ldots,c_{m,\tau_m})      
                                        \right|.
\end{eqnarray*}
But since $(c_{1,\tau_1},\ldots,c_{m,\tau_m})$ is
equidistributed with
$(\theta_1,\ldots,\theta_m)$, the right side of the
displayed inequalities is equidistributed 
with 
$\sup_n \left| 
                                \sum_{m=0}^n 
                                 H_m (d_m (f)) 
                                (\theta_1,\ldots,\theta_m)      
                                        \right|,$
and Step 4 follows.\\
{\bf Proof of Step 3.}  The proof may be done as in
\cite[Theorem~4]{bgs}.  We provide the details to show the
role of analyticity on $\T^N$.  
Here we call a function $\phi\in L^1(\T^N)$
analytic if its Fourier transform is supported in the half-space
$$        \O=\{0\}
\bigcup_{j=1}^N\{(m_1,m_2,\ldots,m_N)\in \Z^N
:\ m_j>0, m_{j+1}=\ldots, m_N=0\}.$$
The following basic properties
of analytic functions on $\T^N$ are easy to prove.
\begin{itemize}
\item
A function $\phi\in L^1(\T^N)$ is analytic if and only
if each term in its martingale difference decomposition,
$d_j(\phi)$ ($j=1,\ldots,N)$,
is analytic in the $j$-th variable $\theta_j$
and has zero mean, i.e.,
$d_j(\phi)\in H^1_0(\T)$.
\item
If $\phi$ is analytic then $\phi^2$ is also analytic.
(This follows from $\O +\O=\O$.)
\item
If $\phi$ is a trigonometric polynomial on $\T^N$, 
then $\phi+i H(\phi)$ is analytic.
\end{itemize}
Getting back to the proof of Step~3, let 
$$g(r_1\theta_1,\dots,r_N\theta_N) = f(r_1\theta_1,\dots,r_N\theta_N) + 
i H (f)(r_1\theta_1,\dots,r_N\theta_N),$$
and let 
$$h = g^2.$$ 
 Both
$g$\ and $h$\ are analytic on $\T^N$.
Hence the
functions $d_m (g)(\theta_1,\ldots , r_m\theta_m)$
and $d_m (h)(\theta_1,\ldots , r_m\theta_m)$
are analytic in the $m$-th variable.
Form the functions $G_{n,t}$\ and $H_{n,t}$\ as 
in (\ref{brownian-martingale}).  By Lemma (\ref{lemma5.2}),
$G_{n,t}$\ and $H_{n,t}$ are martingales relative to 
$\A_{n,t}$. 
We claim that, because of analyticity, we have
\begin{equation}
H_{n,t} = G_{n,t}^2. 
\label{square-of-analytic}
\end{equation} 
To see this, write 
$$ g(\theta_1,\dots,\theta_N) = 
\sum_{k=1}^N d_k(g)(\theta_1,\dots,\theta_k) $$
and
$$ h(\theta_1,\dots,\theta_N) = 
\sum_{k=1}^N d_k(h)(\theta_1,\dots,\theta_k).$$
Then, since all the exponents of $\theta_n$ are positive, 
we get
$$
 \left( \sum_{k=1}^{n-1} d_k(g)(\theta_1,\dots,\theta_k) +
   d_n(g)(\theta_1,\dots,r_n\theta_n) \right)^2
   =
   \sum_{k=1}^{n-1} d_k(h)(\theta_1,\dots,\theta_k) +
   d_n(h)(\theta_1,\dots,r_n\theta_n) $$
and (\ref{square-of-analytic}) easily follows. Consequently,
since the
functions $H_{n,t}$ form a martingale relative
to the $\sigma$-algebra $\A_{n,t}$, we have that
$G_{n,t}^2$ is a martingale relative to this
$\sigma$-algebra.  With this fact in hands,
we can now proceed with
the proof of Step 3 in exactly the same way as in
 \cite[pp. 148-149]{bgs}.  We need a lemma.  
\begin{lemma5.3}
Suppose that $\mu$ and $\nu$ are stopping times with
$\mu\leq \nu$ a.\ e.  Let $f$ be a real-valued
trigonometric polynomial
on $\T^N$
with $\int f dP  = 0$.  Then,
$$\| \tilde{F}_\nu -\tilde{F}_\mu \|_2=
\|  F_\nu -F_\mu \|_2.$$
\label{lemma5.3}
\end{lemma5.3}
{\bf Proof.}  Using the fact that $G_{n,t}^2$
is a martingale, we get
$$
0=\E(G_0^2)=\E(G_\mu^2).$$
Similarly, 
$\E(G_{\nu}^2)=0$.  
Hence, $\E F_\mu^2 = \E \tilde{F}_\mu^2$ and 
$\E F_\nu^2 = \E \tilde{F}_\nu^2$.
Next, we show that 
$\E(F_\mu F_\nu)= \E(F_\mu^2)$, and 
$\E(\tilde{F}_\mu \tilde{F}_\nu)= \E(\tilde{F}_\mu^2)$.
We start with the first equality.  Using Doob's Optional Sampling
Theorem and basic properties of the conditional 
expectation, we see that
$$\E(F_\nu | F_\mu)=F_\mu ,$$ %
$$F_\mu \E(F_\nu | F_\mu)=F_\mu^2 ,$$ %
and so
$$\E(F_\mu F_\nu | F_\mu)=F_\mu^2 .$$ %
Integrating both sides of the
last equality, we get $\E(F_\mu F_\nu)= \E(F_\mu^2)$.
The second equality can be proved similarly.  Thus
\begin{eqnarray*}
\E(F_\mu -F_\nu)^2      &=&     \E F^2_\mu + \E F^2_\nu -2 \E (F_\mu F_\nu)\\
                        &=&     \E F^2_\mu + \E F^2_\nu -2 \E (F_\mu^2) \\
                        &=&     \E F^2_\nu - \E (F_\mu^2)               \\
                        &=&     \E(\tilde{F}_\mu - \tilde{F}_\nu)^2,
\end{eqnarray*}
which completes the proof.\\
The above lemma enables us to establish a 
fundamental inequality.
This is our version of the 
`good $\lambda$' inequality for conjugate functions on $\T^N$.
\begin{lemma5.4}
  With the notation of the previous lemma,
  let $\alpha \ge 1$\ and $\beta > 1$.  Then there is a
constant $c$, depending only on $\alpha$\ and $\beta$, 
such that whenever
$\lambda > 0$\ satisfies
$$ P(G^* > \lambda) \le \alpha P(G^* > \beta \lambda) ,$$
then
$$ P(G^* > \lambda) \le c\,P(c\,F^* > \lambda) .$$
\label{lemma5.4}
\end{lemma5.4}
{\bf Proof.}  
Define stopping times
$$   \mu = \inf\{\, (n,t)\in \cT : |G_{n,t}| > \lambda \} ,\\
   \nu = \inf\{\, (n,t)\in \cT : |G_{n,t}| > \beta \lambda \} . $$
If the set $\{\, (n,t) : |G_{n,t}| > \lambda \}$\ is empty, then we set
$\mu = \infty$.  Otherwise 
$\mu$\ is 
such that $|G_{n,t}| \le \lambda$\
whenever $(n,t) < \mu$, and $|G_\mu| = \lambda$.  
We define $\nu$\ similarly.
Also, we have
that $\mu \le \nu$, that $|G_\mu| = \lambda$\ 
on the set $\{\mu\ne\infty\}
= \{G_\infty^* > \lambda\}$, and that 
$|G_\nu| = \beta\lambda$\ on the set 
$\{\nu\ne\infty\}
= \{G^* > \beta\lambda\}$.  
Thus if $\lambda$\ satisfies the hypothesis of 
the lemma, then
\begin{eqnarray*}
   \E(\chi_{G^*>\lambda} (F_\nu-F_\mu)^2)
  & =  &
         \left\| F_\nu - F_\mu\right\|_2^2 \\
  & =  &
        \frac{1}{2} \left\| G_\nu - G_\mu\right\|_2^2 \\
  & \ge & \frac{1}{2} (\beta \lambda - \lambda)^2 
        P(G^*>\beta \lambda)     \\
  & \ge & c \lambda^2 P(G_\infty^* > \lambda ).  
   \end{eqnarray*}
Also
$$ \E(\chi_{G^*>\lambda} (F_\nu-F_\mu)^4) \le 
   \left\| G_\nu - G_\mu\right\|_4^4
   \le c \lambda^4 P(G_\infty^* > \lambda ). $$
Thus, by a lemma of Paley and Zygmund \cite[Chapter V, (8,26)]{zyg},
$$ P(G^* > \lambda) \le c P(c|F_\nu - F_\mu| > \lambda) .$$
Since $|F_\nu - F_\mu| \le 2 F^*$, the lemma follows.

\bigskip

Now let us finish by proving Step~3.  It is sufficient to show
$\left\| G^* \right\|^*_{1,\infty} \le c \, \left\|
            F^* \right\|^*_{1,\infty}$.  Suppose that 
$$ \left\| G^* \right\|^*_{1,\infty} = 
   \sup_{\lambda>0} \lambda P(G^*>\lambda)
   = A .$$
Pick $\lambda_0$\ such that 
$2\lambda_0 P(G^*>2\lambda_0) \ge A/2$.  
Then
$\lambda_0 P(G^*>\lambda_0) \le A$, 
and thus $\lambda_0$\ satisfies
the hypothesis of the lemma with 
$\alpha = 4$\ and $\beta = 2$.  Then
it follows that 
$$ \| F^*\|^*_{1,\infty} \ge 
\lambda_0 P(c F^* > \lambda_0) \ge c A/4 ,$$
as desired.

{\bf Acknowledgements}  The research of the authors was supported 
by grants from the National Science Foundation (U.\ S.\ A.).

\end{document}